%% file: 2013_hygroplate.tex
\documentclass{TTP_DSL2006}

\usepackage{graphicx,graphics}
\usepackage[intlimits]{amsmath}
\usepackage{amssymb}
\usepackage{exscale}
\usepackage{times}
\usepackage{subfigure}
\usepackage{subfig}
\usepackage{multirow}
\usepackage[abs]{overpic}

\usepackage{pgfplots}
\usetikzlibrary{plotmarks}

\renewcommand{\large}{\fontsize{14}{18pt}\selectfont}
\renewcommand{\small}{\fontsize{11}{13.6pt}\selectfont}

\newcommand{\titleformat}{\sffamily\bfseries \large}						
\newcommand{\authorformat}{\sffamily \large}							
\newcommand{\keywordsformat}{\noindent \small \sffamily}				
\newcommand{\abstractformat}{\noindent \textbf}						
\newcommand{\contentformat}{\rmfamily \normalsize\vspace{18pt}}			
\renewcommand{\subsection}{\textbf}	

\newcommand{\Eref}[1]{Equation (\ref{#1})}
\newcommand{\fref}[1]{Figure (\ref{#1})}
\newcommand{\Erefs}[1]{Equations (\ref{#1})}

\newcommand{\xx}{\mathbf{x}}
\newcommand{\KK}{\mathbf{K}}
\newcommand{\bm}{\mathbf{M}}

\newcommand{\ff}{\mathbf{f}}

\newcommand{\rmd}{\rm{d}}

\newcommand{\bveps}{\boldsymbol{\varepsilon}}
\newcommand{\bvsig}{\boldsymbol{\sigma}}

\begin{document}

\title{\titleformat Hygrothermal effects on free vibration and buckling of laminated composites with cutouts}

\author{\authorformat Sundararajan Natarajan \inst{1}$^{,\rm{a}}$\text{,} and Pratik S Deogekar \inst{2}$^{,\rm{b}}$}

\institute{\sffamily School of Civil \& Environmental Engineering, The University of New South Wales, Sydney, NSW 2052, Australia. \and Department of Civil Engineering, Indian Institute of Technology, Bombay.}



\maketitle


\keywordsformat{{\textbf{Keywords:}}  Vibration, Buckling, Reissner-Mindlin plate, extended finite element method, hygrothermal effects.}

\contentformat

\abstractformat{Abstract.} The effect of moisture concentration and the thermal gradient on the free flexural vibration and buckling of laminated composite plates are investigated. The effect of a centrally located cutout on the global response is also studied. The analysis is carried out within the framework of the extended finite element method. A Heaviside function is used to capture the jump in the displacement and an enriched shear flexible 4-noded quadrilateral element is used for the spatial discretization. The formulation takes into account the transverse shear deformation and accounts for the lamina material properties at elevated moisture concentrations and temperature. The influence of the plate geometry, the geometry of the cutout, the moisture concentration, the thermal gradient and the boundary conditions on the free flexural vibration is numerically studied. 

\section{Introduction}
\vspace{-6pt}
Fibre reinforced laminated composites belong to a class of engineered materials that has found increased utilization as structural elements in the construction of aeronautical and aerospace vechicles, sports, as well as civil and mechanical structures. This is because of their excellent strength-to and stiffness-to-weight ratios and the possibility to tailor their properties to optimize the structural response. However, the analysis of such structures is very demanding due to coupling between membrane, torsion and bending strains; weak transverse rigidities; and discontinuity of the mechanical characteristics through the thickness of the laminates. The application of analytical/numerical methods based on various 2D theories have attracted the attention of the research community. In general, three different approaches have been used to study laminated composite structures: single layer theories, discrete layer theories and mixed plate theory. In the single layer theory approach, layers in laminated composites are assumed to be one equivalent single layer (ESL), whereas in the discrete layer theory approach, each layer is considered in the analysis. Although the discrete layer theories provide very accurate prediction of the displacements and the stresses, increasing the number of layers increases the number of unknowns. This can be prohibitively costly and significantly increase the computational time~\cite{wuchen2005}. To overcome the above limitation, zig-zag models developed by Murukami~\cite{murukami1986} can satisfy the transverse shear stresses continuity conditions at the interfaces. Moreover, the number of unknowns are independent of the number of layers. Carrera~\cite{carrerademasi2002,carrerademasi2002a,carrera2003} derived a series of axiomatic approaches, coined as `Carrera Unified Formulation' (CUF) for the general description of two-dimensional formulations for multilayered plates and shells. With this unified formulation it is possible to implement in a single software a series of hierarchical formulations, thus affording a systematic assessment of different theories, ranging from simple ESL models up to higher order layerwise descriptions. This formulation is a valuable tool for gaining a deep insight into the complex mechanics of laminated structures.

Plates with cutouts are extensively used in transport vehicle structures. Cutouts are made to lighten the structure, for ventilation, to provide accessibility to other parts of the structures and for altering the resonant frequency. Therefore, the natural frequencies of plates with cutouts are of considerable interest to designers of such structures. Most of the earlier investigations on plates with cutouts have been confined to isotropic plates~\cite{paramasivam1973,aliatwal1980,huangsakiyama1999} and laminated composites~\cite{reddy1982,sivakumariyengar1998}. Moreover, the laminated composites may be subjected to moisture and temperature environment during its service life. The moisture concentration and thermal environment can have significant impact on the response of such laminated structures. Whitney and Ashton~\cite{whitneyashton1971} employed Ritz method to analyze the effect of environment on the free vibration of symmetric laminates. Patel \textit{et al.,}~\cite{patelganapathi2002} employed shear flexible Q8 quadrilateral element to study the hygrothermal effects on the structural behaviour of thick composite laminates. Patel \textit{et al.,} employed higher order accurate theory and studied the importance of retaining higher order terms in the displacement approximation. The effect of cutouts on the buckling behaviour of laminated composites were studied in~\cite{nemeth1988,sairamsinha1992}. And more recently, Komur \textit{et al.,}~\cite{komursen2010} and Ghannadpour \textit{et al.,}~\cite{ghannadpournajafi2006} studied the buckling behaviour of laminated composites with circular and elliptical cutouts using the finite element method and first order shear deformation theory. Their study was restricted to a limited number of configurations, because the mesh has to conform to the geometry. Moreover, to the author's knowledge the effect of cutout on the free vibration and buckling behaviour of laminated composites in hygrothermal environment has not been studied earlier or was limited to simple configurations. In this study, we present a framework that provides flexibility to handle internal discontinuities. 

In this paper, we study the influence of a centrally located cutout on the fundamental natural frequency and the critical load of multilayered composite laminated plates in hygrothermal environment. Circular and elliptical cutouts are considered for the study. A structured quadrilateral mesh is used and the cutouts are modelled independent of the mesh within the extended finite element (XFEM) framework. A systematic parametric study is carried out to bring the effect of the boundary conditions, the thermal gradient $\Delta T$, the change in moisture concentration $\Delta C$, the geometry of the cutout on the free flexural vibration and buckling of laminated composites. \\


\section{Theoretical formulation}
\label{theory}
\input{platetheory}

\section{Spatial discretization}
\label{spatdisc}
\input{spatialdis}

\section{Numerical Examples}
\label{numerics}

\input{Validation}

\input{vibration}

\input{buck}

\section{Conclusion}
The extended finite element framework was adopted to study the hygrothermal effects on the free vibration and buckling of multilayered laminated composites with a centrally located cutout. The formulation developed is general in nature and can handle non-uniform distributions of moisture and temperature, although only uniform distribution is considered in the present analysis. The broad conclusion that can be made from this parametric study is that with the increase in the uniform moisture concentration and the temperature, the reduction in the fundamental natural frequency and the critical load need not be linear and could lead to instability depending on the value of the moisture content, temperature and side-to-thickness ratio and aspect ratio. It can also be concluded that the presence of moisture content has negligible effect on the fundamental frequency of thick laminated plate.

\newpage
\bibliographystyle{plain}
\bibliography{myRefFGM}

\end{document}

%% file: platetheory.tex
\begin{figure}[htpb]
\centering
\includegraphics[scale=0.75]{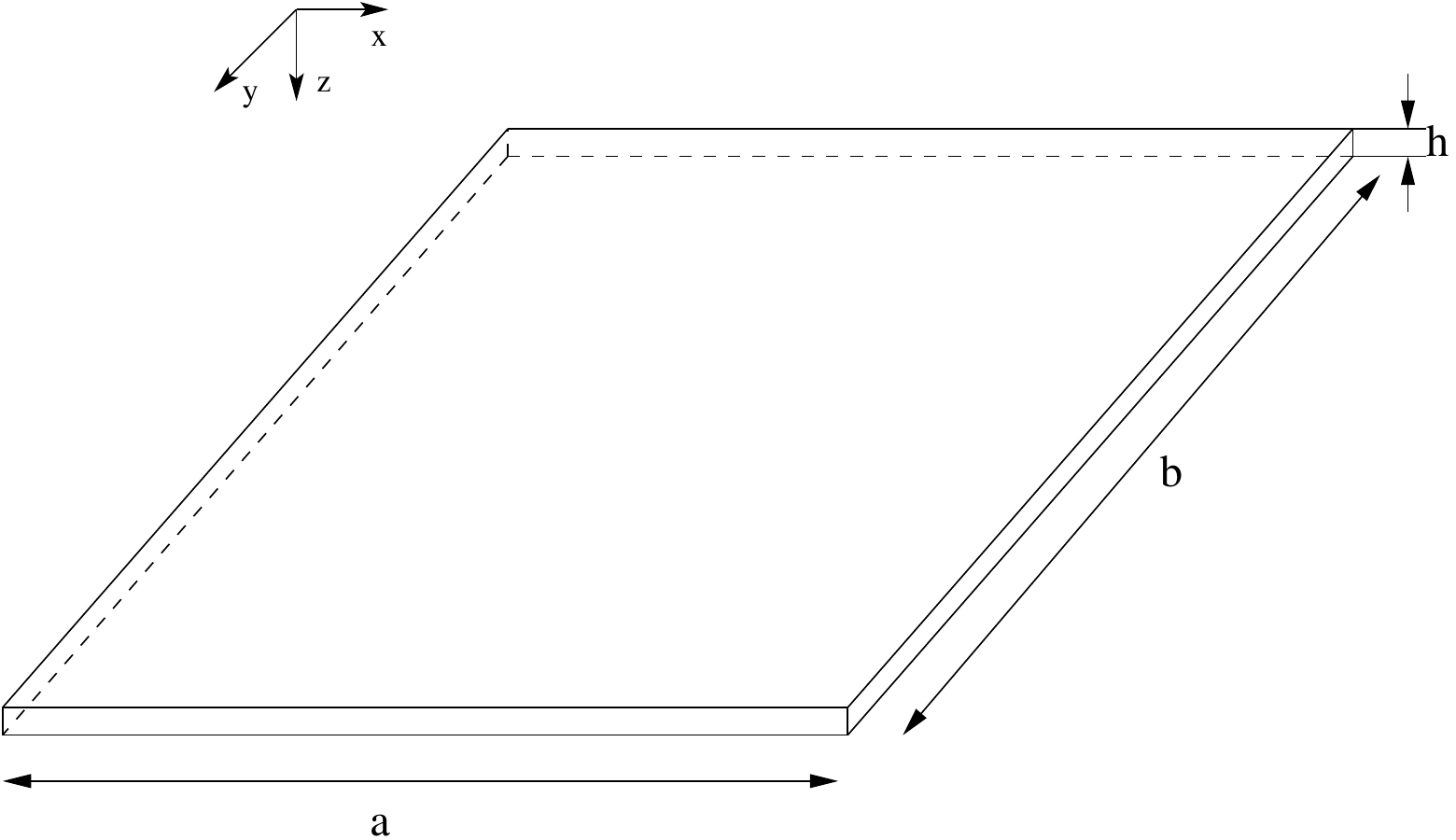}
\caption{Coordinate system of a rectangular laminated plate.}
\label{fig:platefig}
\end{figure}

Using the Mindlin formulation, the displacements $u,v,w$ at a point $(x,y,z)$ in the plate (see \fref{fig:platefig}) from the medium surface are expressed as functions of the mid-plane displacements $u_o,v_o,w_o$ and independent rotations $\beta_x,\beta_y$ of the normal in $yz$ and $xz$ planes, respectively, as
\begin{eqnarray}
u(x,y,z,t) &=& u_o(x,y,t) + z \beta_x(x,y,t) \nonumber \\
v(x,y,z,t) &=& v_o(x,y,t) + z \beta_y(x,y,t) \nonumber \\
w(x,y,z,t) &=& w_o(x,y,t) 
\label{eqn:displacements}
\end{eqnarray}
where $t$ is the time. The strains in terms of mid-plane deformation can be written as:
\begin{equation}
\bveps  = \left\{ \begin{array}{c} \bveps_p \\ \mathbf{0} \end{array} \right \}  + \left\{ \begin{array}{c} z \bveps_b \\ \bveps_s \end{array} \right\} - \left\{ \overline{\bveps}_o \right\}
\label{eqn:strain1}
\end{equation}

The midplane strains $\bveps_p$, the bending strains $\bveps_b$ and the shear strains $\varepsilon_s$ in \Eref{eqn:strain1} are written as:
\begin{equation}
\renewcommand{\arraystretch}{1}
\bveps_p = \left\{ \begin{array}{c} u_{o,x} \\ v_{o,y} \\ u_{o,y}+v_{o,x} \end{array} \right\}, \hspace{0.2cm}
\renewcommand{\arraystretch}{1}
\bveps_b = \left\{ \begin{array}{c} \beta_{x,x} \\ \beta_{y,y} \\ \beta_{x,y}+\beta_{y,x} \end{array} \right\}, \hspace{0.2cm}
\renewcommand{\arraystretch}{1}
\bveps_s = \left\{ \begin{array}{c} \beta _x + w_{o,x} \\ \beta _y + w_{o,y} \end{array} \right\}.
\label{eqn:platestrain}
\end{equation}
where the subscript `comma' represents the partial derivative with respect to the spatial coordinate succeeding it. The strain vector $\left\{ \overline{\bveps}_o \right\}$ due to temperature and moisture is represented as:

\begin{equation}
\overline{\bveps}_o = \left\{ \begin{array}{c} \overline{\varepsilon}_{xx} \\ \overline{\varepsilon}_{yy} \\ \overline{\varepsilon}_{xy} \end{array} \right\} = \Delta T \left\{ \begin{array}{c} \alpha_x \\ \alpha_y \\ \alpha_{xy} \end{array} \right\} + \Delta C \left\{ \begin{array}{c} \beta_x \\ \beta_y \\ \beta_{xy} \end{array} \right\}
\end{equation}
where $\Delta T$ and $\Delta C$ are the rise in temperature and the moisture concentration, respectively. $\alpha_x, \alpha_y$ and $\alpha_{xy}$ are the thermal expansion coefficients in the plate coordinates and can be related to the thermal coefficients $(\alpha_1, \alpha_2 \textup{and} \alpha_3$ in the material principal directions and $\beta_x, \beta_y$ and $\beta_{xy}$ are the moisture expansion coefficients similar to thermal expansion coefficients in the plate coordinates. The constitutive relations for an arbitrary layer $k$ in the laminate $(x,y,z)$ coordinate system can be expressed as:

\begin{equation}
\bvsig = \left\{ \begin{array}{c} \sigma_{xx} \\ \sigma_{yy} \\ \sigma_{xy} \end{array} \right\} = [\overline{Q}_k ] \left\{ \left\{ \begin{array}{c} \bveps_p \\ \mathbf{0} \end{array} \right \}  + \left\{ \begin{array}{c} z \bveps_b \\ \bveps_s \end{array} \right\} - \left\{ \overline{\bveps}_o \right\} \right\}
\end{equation}
where the terms of $[ \overline{Q}_k]$ matrix of $k^{\rm th}$ ply are referred to the laminate axes and can be obtained from the $[ Q_k]$ corresponding to the fibre directions with the appropriate transformations. The governing equations are obtained by applying Lagrangian equations of motion:
\begin{equation}
\frac{d}{dt} \left[ \frac{\partial(T-U)}{\partial \dot{\boldsymbol{\delta}}_i} \right] - \left[ \frac{\partial(T-U)}{\partial \boldsymbol{\delta}_i} \right] = 0, \hspace{0.2cm} i=1, \cdots n
\end{equation}
where $T$ is the kinetic energy, given by:
\begin{equation}
T(\boldsymbol{\delta}) = {1 \over 2} \int_{\Omega} \left\{p (\dot{u}_o^2 + \dot{v}_o^2 + \dot{w}_o^2) + I(\dot{\theta}_x^2 + \dot{\theta}_y^2) \right\}~\rmd \Omega
\label{eqn:kinetic}
\end{equation}
where $p = \int_{-h/2}^{h/2} \rho~\rmd z$, $I = \int_{-h/2}^{h/2} z^2 \rho~\rmd z$ and $\rho$ is the mass density of the plate. The strain energy function $U$ is given by:
\begin{equation}
U(\boldsymbol{\delta}) = \frac{1}{2} \iint \left[ \sum\limits_{k=1}^n \int\limits_{h_k}^{h_{k+1}} \bvsig^{\rm T} \bveps~\rmd z \right] \rmd x \rmd y
\label{eqn:potential}
\end{equation}
where $\boldsymbol{\delta} = \{u,v,w,\beta_x,\beta_y\}$ is the vector of the degrees of freedom associated to the displacement field in a finite element discretization. Substituting \Erefs{eqn:kinetic} - (\ref{eqn:potential}) in Lagrange's equations of motion and following the procedure given in~\cite{rajasekaran1973}, the following discretized equation is obtained:

\begin{equation}
\bm \ddot{\boldsymbol{\delta}} + \left[ \KK + \KK_R + \KK_G \right]\boldsymbol{\delta}  = \ff_T
\end{equation}
where $\KK$ is the global linear stiffness matrix, $\KK_R$ and $\KK_G$ are the global geometric stiffness due to the residual stresses and the applied in-plane mechanical loads, respectively, $\bm$ is the global mass matrix and $\ff_T$ is the global hygrothermal load vector. After substituting the characteristic of the time function~\cite{Ganapathi1991} $\ddot{\boldsymbol{\delta}} = -\omega^2 \boldsymbol{\delta}$, the following algebraic equation is obtained:

\paragraph{Static bending:} $ \KK \boldsymbol{\delta} = \ff_T$

\paragraph{Free vibration:} $\left[ \left( \KK + \KK_R \right) - \omega^2 \bm \right] \boldsymbol{\delta} = \mathbf{0}$

\paragraph{Buckling:} $\left[ \left( \KK + \KK_R \right) - \lambda \KK_G \right] \boldsymbol{\delta} = \mathbf{0}$

where $\omega$ is the natural frequency and $\lambda$ is the buckling load. The residual stress state depends on the ply lay-up. Hence, to evaluate the stress state, pre-buckling displacement field for the assumed hygro-thermal-mechanical load is obtained by solving static bending. The displacement field is then used to calculate the stresses and in turn, $\KK_R$ and $\KK_G$ in matrices.

%% file: spatialdis.tex
The plate element employed here is a $\mathcal{C}^o$ continuous shear flexible field consistent element with five degrees of freedom $(u_o,v_o,w_o,\beta_x,\beta_y)$ at four nodes in a 4-noded quadrilateral (QUAD-4) element. The displacement field within the element is approximated by:
\begin{equation}
\{ u_o^e,v_o^e,w_o^e,\beta_x^e,\beta_y^e\} = \sum\limits_{J=1}^4 N_J \{u_{oJ}, v_{oJ}, w_{oJ},\beta_{xJ},\beta_{yJ} \},
\end{equation} 
where $u_{oJ}, v_{oJ}, w_{oJ},\beta_{xJ},\beta_{yJ}$ are the nodal variables and $N_J$ are the shape functions for the bi-linear QUAD-4 element. If the interpolation functions for a QUAD-4 are used directly to interpolate the five variables $(u_o,v_o,w_o,\beta_x,\beta_y)$ in deriving the shear strains and the membrane strains, the element will lock and show oscillations in the shear and the membrane stresses. The oscillations are due to the fact that the derivative functions of the out-of plate displacement, $w_o$ do not match that of the rotations ($\beta_x, \beta_y$) in the shear strain definition, given by \Eref{eqn:platestrain}. To alleviate the locking phenomenon, the terms corresponding to the derivative of the out-of plate displacement, $w_o$ must be consistent with the rotation terms, $\beta_x$ and $\beta_y$. The present formulation, when applied to thin plates, exhibits shear locking. In this study, field redistributed shape functions are used to alleviate the shear locking.~\cite{Somashekar1987,natarajanbaiz2011a} The field consistency requires that the transverse shear strains and the membrane strains must be interpolated in a consistent manner. Thus, the $\beta_x$ and $\beta_y$ terms in the expressions for the shear strain $\bveps_s$ have to be consistent with the derivative of the field functions, $w_{o,x}$ and $w_{o,y}$. 

\subsection{Representation of discontinuity surface}
The finite element framework requires the underlying finite element mesh to conform to the discontinuity surface. The recent introduction of implicit boundary definition-based methods, viz., the extended/generalized FEM (XFEM/GFEM), alleviates the shortcomings associated with the meshing of the discontinuity surface. In this study, the partition of unity framework is employed to represent the discontinuity surface independent of the underlying mesh.

\begin{equation}
(u^h,v^h,w^h,\beta_x^h,\beta_y^h)\left(\xx\right) = \underbrace{ \sum_{I \in \mathcal{N}^{\rm{fem}}} N_I(\xx) (u_I^s,v_I^s,w_I^s,\beta_{xI}^s,\beta_{yI}^s)}_{\rm FEM} + \underbrace{ \sum_{J \in
\mathcal{N}^{\rm{c}}} N_J(\xx) H(\xx) (b_J^u,b_J^v,b_J^w, b_J^{\beta_x},b_J^{\beta_y})}_{\rm Enriched~part} 
\label{eqn:platexfem1}
\end{equation}
where $\mathcal{N}^{\rm{fem}}$ is a set of all the nodes in the finite element mesh and $\mathcal{N}^{\rm{c}}$ is a set of nodes that are enriched with the Heaviside function. In \Eref{eqn:platexfem1}, $(u_I^s,v_I^s,w_I^s,\beta_{x_I}^s,\beta_{y_I}^s)$ are the nodal unknown vectors associated with the continuous part of the finite element solution, $b_J$ is the nodal enriched degree of freedom vector associated with the Heaviside (discontinuous) function. In this study, a level set approach is followed to model the cutouts. The geometric interface (for example, the boundary of the cutout) is represented by the zero level curve $\phi \equiv \phi(\xx,t) = 0$. The interface is located from the value of the level set information stored at the nodes. The standard FE shape functions can be used to interpolate $\phi$ at any point $\xx$ in the domain as:
\begin{equation}
\phi(\xx) = \sum\limits_I N_I(\xx) \phi_I
\end{equation}
where the summation is over all the nodes in the connectivity of the elements that contact $\xx$ and $\phi_I$ are the nodal values of the level set function. For circular cutout, the level set function is given by:
\begin{equation}
\phi_I = || \xx_I - \xx_c|| - r_c
\end{equation}
where $\xx_c$ and $r_c$ are the center and the radius of the cutout. For an elliptical cutout oriented at an angle $\theta$, measured from the $x-$ axis the level set function is given by:
\begin{equation}
\phi_I = \sqrt{a_1(x_I-x_c)^2-a_2(x_I-x_c)(y_I-y_c) + a_3(y_I-y_c)^2} - 1 
\end{equation}
where
\begin{equation}
a_1 = \left( {\cos \theta \over d} \right)^2, \hspace{0.15cm} a_2 = 2\cos\theta \sin\theta\left( {1 \over d^2} - {1 \over e^2} \right), \hspace{0.15cm} a_3 = \left( {\sin \theta \over d} \right)^2 + \left( {\cos \theta \over e}\right)^2.
\end{equation}
where $d$ and $e$ are the major and minor axes of the ellipse and $(x_c,y_c)$ is the center of the ellipse. 

\subsection{Numerical integration over enriched elements}
A consequence of adding custom tailored enrichment functions to the FE approximation basis, which are not necessarily smooth functions is that, special care has to be taken in numerically integrating over the elements that are intersected by the discontinuity surface. The standard Gau\ss~quadrature cannot be applied in elements enriched by discontinuous terms, because Gau\ss~quadrature implicitly assumes a polynomial approximation. One potential solution for the purpose of numerical integration is by partitioning the elements into subcells (to triangles for example) aligned to the discontinuous surface in which the integrands are continuous and differentiable~\cite{Belytschko2009}. The other techniques that can be employed are Schwarz Christoffel Mapping~\cite{natarajanbordas2009,natarajanmahapatra2010}, Generalized quadrature~\cite{mousavisukumar2011} and Smoothed eXtended FEM~\cite{bordasnatarajan2011}. In the present study, a triangular quadrature with sub-division is employed along with the integration rules described in Table \ref{table:subcellgausspt}. For the elements that are not enriched, a standard 2 $\times$ 2 Gaussian quadrature rule is used.

\begin{table}[htpb]
\caption{Integration rules for enriched and non-enriched elements in the presence of a crack}
\centering
\begin{tabular}{lr}
\hline
Element Type & Gau\ss ~points\\
\hline
Non-enriched element & 4  \\
Tip element & 13 per triangle \\
Tip blending element & 16 \\
Split element & 3 per triangle  \\
Split blending element & 4 \\
Split-Tip blending element & 4 per triangle \\
\hline
\end{tabular}
\label{table:subcellgausspt}
\end{table}

%% file: Validation.tex
\begin{figure}[htpb]
\centering
\includegraphics[scale=0.6]{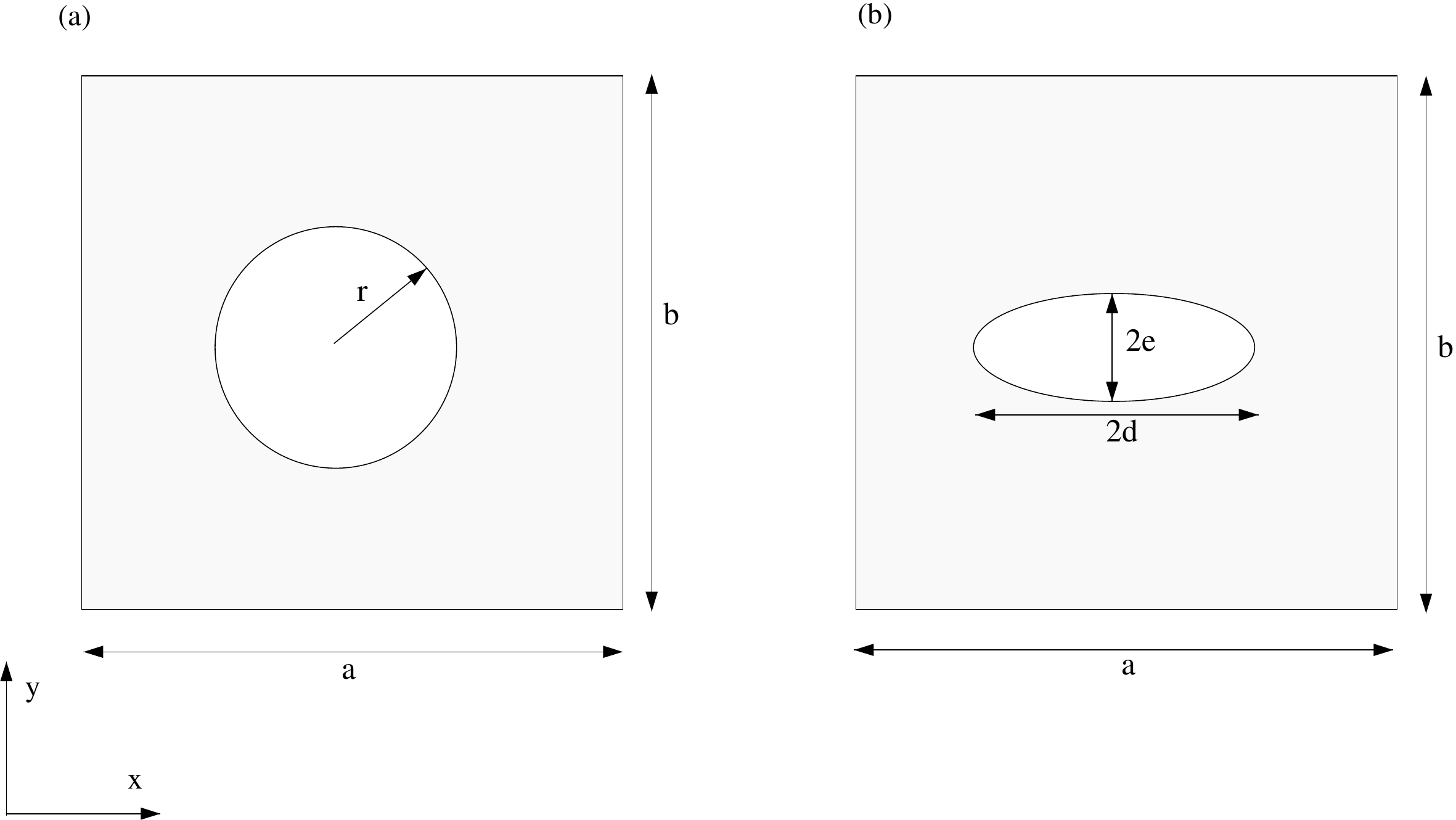}
\caption{Plate with a centrally located circular and an elliptical cutout. $r$ is the radius of the circular cutout, $2d$ and $2e$ are the major and minor axes defining the ellipse.}
\label{fig:problemdescription}
\end{figure}

In this section, we study the influence of a centrally located cutout on the fundamental frequencies of FGM plates. We consider both square and rectangular plates with simply supported and clamped boundary conditions. Two different cutout shapes, viz., circular and elliptical cutouts (see \fref{fig:problemdescription}) are considered in this study. Although the formulation presented here is general, the analysis is carried out for cross-ply laminates subjected to uniform distributions of moisture and/or temperature. The lamina properties at the elevated moisture concentration and temperature is given in Tables \ref{table:moistproperty} and \ref{table:tempproperty}. The effect of the plate slenderness ratio $a/h$, the  plate aspect ratio $b/a$, the cutout radius $r/a$, the cutout geometry $d/e$ and the boundary condition on the natural frequencies and the critical buckling load are numerically studied. The boundary conditions for simply supported and clamped cases are :

\noindent \emph{Simply supported boundary condition}: \\
\begin{equation}
u_o = w_o = \theta_y = 0 \hspace{0.2cm} ~\textup{on} \hspace{0.2cm}  x=0,a; \hspace{0.2cm}
v_o = w_o = \theta_x = 0 \hspace{0.2cm} ~\textup{on} \hspace{0.2cm}  y=0,b
\end{equation}

\noindent \emph{Clamped boundary condition}: \\
\begin{equation}
u_o = w_o = \theta_y = v_o = \theta_x = 0 \hspace{1cm} ~\textup{on} ~ x=0,a \hspace{0.2cm} \& \hspace{0.2cm}  y=0,b
\end{equation}

\begin{table}
\renewcommand\arraystretch{1}
\caption{Elastic moduli of graphite/epoxy lamina at different moisture concentrations, $G_{13}=G_{12}, G_{23} = 0.5G_{12}, \nu_{12}=0.3, \beta_1 = 0$ and $\beta_2 = 0.44$.}
\centering
\begin{tabular}{lrrrrrrr}
\hline
Elastic & \multicolumn{7}{l}{Moisture concentration C($\%)$} \\
\cline{2-8}
Moduli (GPa) & 0.0 & 0.25 & 0.50 & 0.75 & 1.00 & 1.25 & 1.50 \\
\hline
$E_1$ & 130 & 130 & 130 & 130 & 130 & 130 & 130 \\
$E_2$ & 9.50 & 9.25 & 9.00 & 8.75 & 8.50 & 8.50 & 8.50 \\
$G_{12}$ & 6.0 & 6.0 & 6.0 & 6.0 & 6.0 & 6.0 & 6.0\\
\hline
\end{tabular}
\label{table:moistproperty}
\end{table}

\begin{table}
\renewcommand\arraystretch{1}
\caption{Elastic moduli of graphite/epoxy lamina at different temperatures, $G_{13}=G_{12}, G_{23} = 0.5G_{12}, \nu_{12}=0.3, \alpha_1 = -0.3\times$10$^{-6}$/K  and $\alpha_2 = 28.1\times$10$^{-6}$/K.}
\centering
\begin{tabular}{lrrrrrr}
\hline
Elastic & \multicolumn{6}{l}{Temperature $T$ (K)} \\
\cline{2-7}
Moduli (GPa) & 300 & 325 & 350 & 375 & 400 & 425 \\
\hline
$E_1$ & 130 & 130 & 130 & 130 & 130 & 130 \\
$E_2$ & 9.50 & 8.50 & 8.00 & 7.50 & 7.00 & 6.75 \\
$G_{12}$ & 6.0 & 6.0 & 5.5 & 5.0 & 4.75 & 4.50 \\
\hline
\end{tabular}
\label{table:tempproperty}
\end{table}

\paragraph{Validation} Before proceeding with a detailed study on the effect of different parameters on the natural frequency and the critical load, the formulation developed herein is validated against available closed form/analytical solutions. The critical load and the natural frequency of a cross-ply laminate exposed to moisture and temperature are presented in Table \ref{table:freqbuckvalidation}, along with the Ritz solutions~\cite{whitneyashton1971} and with Q8 element~\cite{patelganapathi2002}. It can be seen that the results from the present formulation compare very well with the available solutions and based on a progressive mesh refinement, a structured quadrilateral mesh of 40$\times$ 40 is found to be adequate to model the full laminate for the present analysis. From the Table \ref{}, it is seen that the percentage difference between a structured 30 $\times$ 30 and 40 $\times$ 40 mesh is less than 0.1\%, hence, a structured mesh of 30 $\times$ 30 is used for the analysis. Next, through the present formulation, the influence of various parameters on the natural frequency and the critical load is studied.
\begin{table}
\renewcommand\arraystretch{1}
\caption{Comparison of natural frequency and critical load for four layered cross-ply laminates with $a/h=$ 100.}
\centering
\begin{tabular}{lrrrrrrr}
\hline
Mesh & \multicolumn{2}{c}{Frequency, $\Omega=\omega a^2 \sqrt{\frac{\rho}{(E_2h^2)}}$} && \multicolumn{2}{c}{ $\overline{N}_{xx}^{\rm cr} = N_{xx}^{\rm cr}/N_{xx}^{\ast {\rm cr}}$} \\
\cline{2-3} \cline{5-6}
 & $C=$ 0.1$\%$ & $T=$ 325K && $C=$ 0.1$\%$ & $T=$ 325K \\
\hline
10 $\times $ 10 & 9.6133 & 8.2604 && 0.6158 & 0.4571 \\
20 $\times$ 20 & 9.4596 & 8.0926 && 0.6100 & 0.4488 \\
30 $\times$ 30 & 9.4345 & 8.0651 & & 0.6090 & 0.4475 \\
40 $\times$ 40 & 9.4260 & 8.0559 && 0.6087 & 0.4393 \\
Ref.~\cite{whitneyashton1971} & 9.4110 & 8.0680 && 0.6091 & 0.4477\\
Ref.~\cite{patelganapathi2002} & 9.3993 & 8.0531 && 0.6084 & 0.4466\\
\hline
\end{tabular}
\label{table:freqbuckvalidation}
\end{table}

%% file: vibration.tex
\paragraph{Vibration}
Consider a plate with side lengths $a$ and $b$ and thickness $h$. A laminated plate with ply sequence (0$^\circ$/90$^\circ$/90$^\circ$/0$^\circ$) is considered for the analysis. In all cases, we present the non dimensionalized free flexural frequencies as, unless specified otherwise:
\begin{equation}
\Omega = \omega \left( \frac{a^2}{h} \right) \sqrt{ \frac{\rho}{E_2}}
\label{eqn:nondimfreq}
\end{equation}

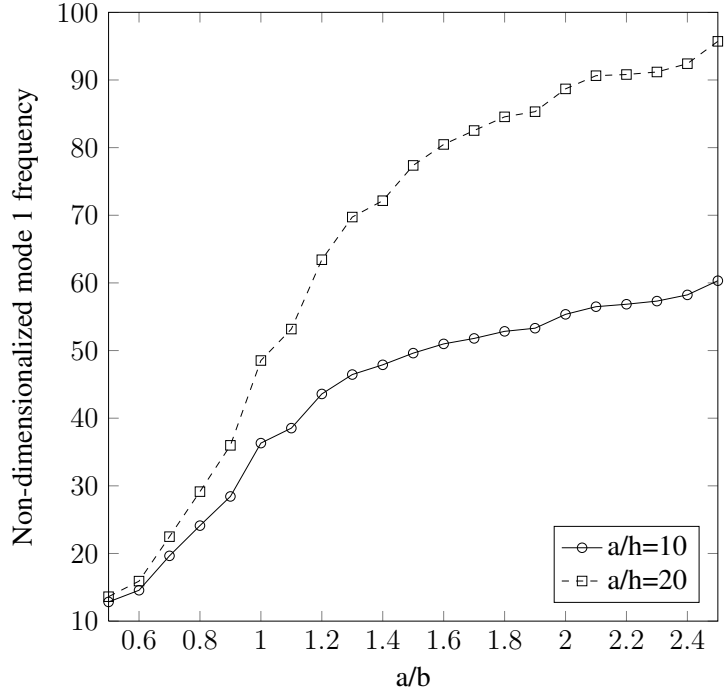
\begin{figure}
\centering
\newlength\figureheight 
\newlength\figurewidth 
\setlength\figureheight{9cm} 
\setlength\figurewidth{9cm}
\input{abVsfreq.tikz}
\caption{Influence of the plate aspect ratio on the fundamental frequency of a simply supported laminated plate with (0$^\circ$/90$^\circ$/90$^\circ$/0$^\circ$) and a centrally located circular cutout, $r_o/a=$ 0.2. The laminated plate is exposed to a moisture concentration, $C=$ 1\% and $T=$ 300K.}
\label{fig:abratio}
\end{figure}

\begin{figure}
\centering
\setlength\figureheight{8cm} 
\setlength\figurewidth{8cm}
\input{ravsfreq.tikz}
\caption{Non-dimensionalized mode 1 frequency as a function of cutout radius $r_o/a$ for a simply supported square laminated plate with $a/h=$10 and (0$^\circ$/90$^\circ$/90$^\circ$/0$^\circ$). The plate is exposed to different thermal environment and moisture concentration $C=$ 0\%.}
\label{fig:tempeffect}
\end{figure}
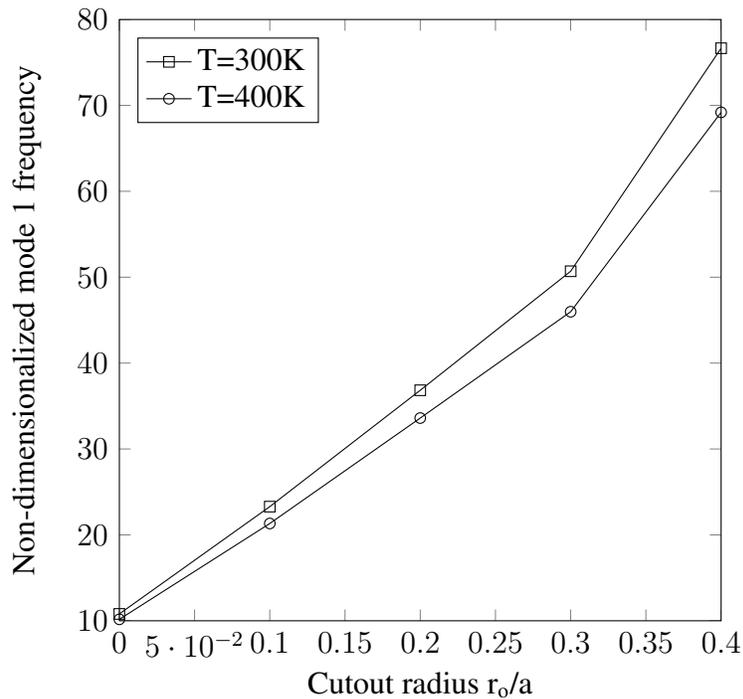

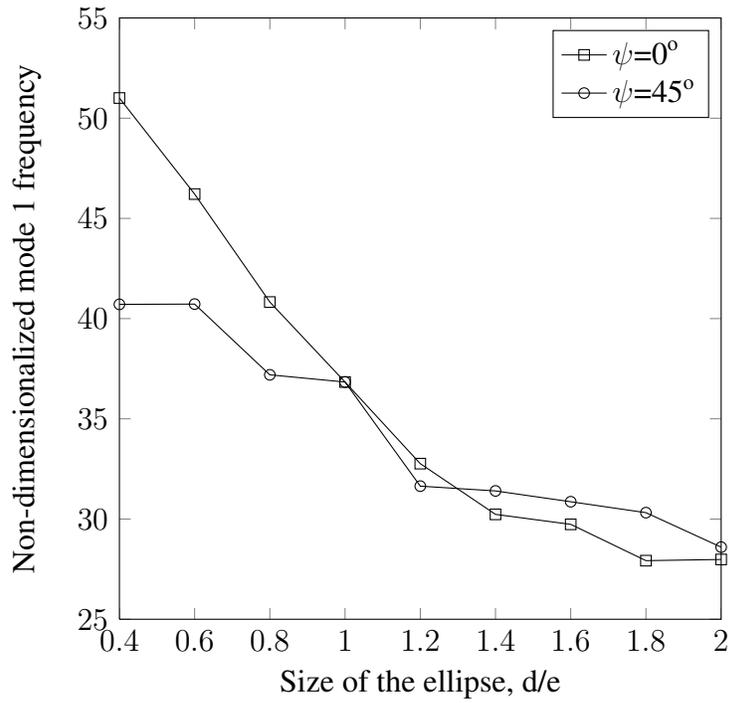
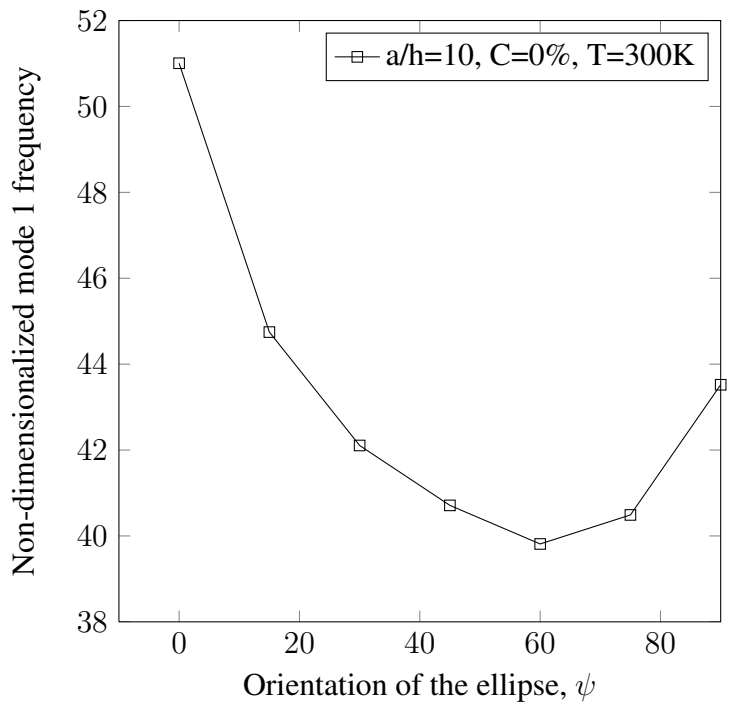
\begin{figure}
\centering
\setlength\figureheight{8cm} 
\setlength\figurewidth{8cm}
\subfigure[Influence of the geometry of the cutout $d/e$]{\input{devsfreq.tikz}}
\subfigure[Influence of the orientation of the cutout $\psi$]{\input{deOrienvsfreq.tikz}}
\caption{Influence of the geometry and the orientation of the cutout on the linear frequency for a square (0$^\circ$/90$^\circ$/90$^\circ$/0$^\circ$) plate with $a/h=$10. The plate is exposed to temperature $T=$ 300K and mositure concentration, $C=$ 0\%.}
\label{fig:elliporientation}
\end{figure}

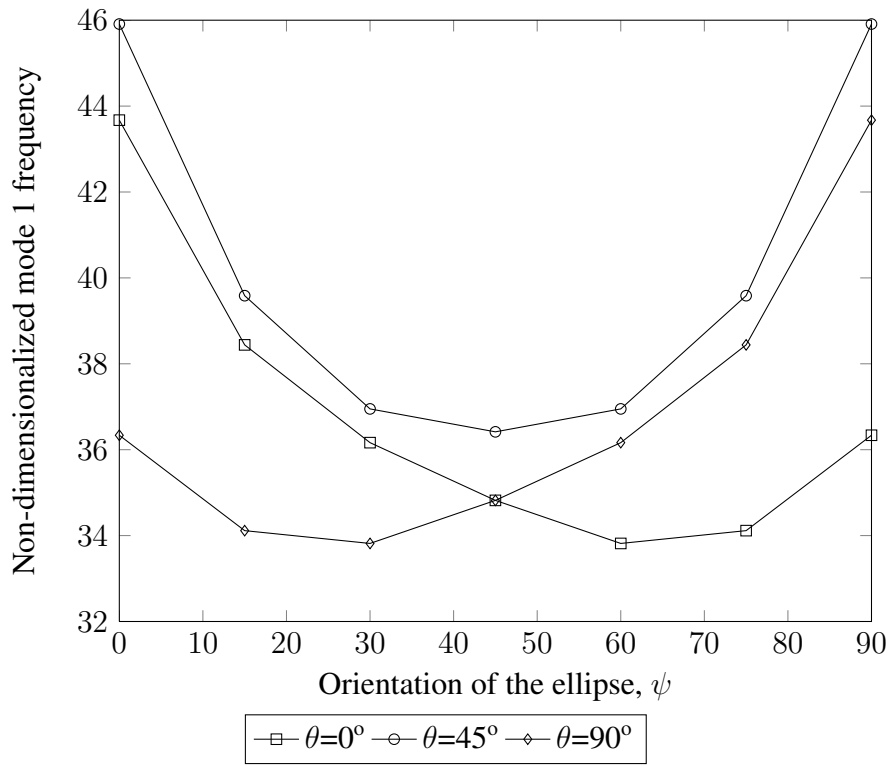
\begin{figure}
\centering
\setlength\figureheight{8cm} 
\setlength\figurewidth{10cm}
\input{singlevsfreq.tikz}
\caption{Non-dimensionalized mode 1 frequency as a function of orientation of the cutout $\psi$ for a simply supported square laminated plate with $a/h=$10 and for various ply orientations. In this case, only one single lamina is considered and the plate is exposed to temperature $T=$ 300K and moisture concentration $C=$ 0\%.}
\label{fig:plyorientation}
\end{figure}

\fref{fig:abratio} shows the influence of plate aspect ratio $a/b$ on the natural frequency for a simply supported laminated plate (0$^\circ$/90$^\circ$/90$^\circ$/0$^\circ$) with a centrally located circular cutout $r_o/a=$ 0.2 and exposed to moisture concentration $C=$1\% and temperature $T=$ 300K. The effect of plate thickness is also shown in \fref{fig:abratio}. It is seen that increasing the plate aspect ratio and decreasing the thickness of the plate, increases the natural frequency. The effect of moisture concentration $C$, the plate thickness $a/h$ and the boundary conditions is shown in \fref{fig:moistureeffect}.  Increasing the moisture concentration has greater impact for a plate with smaller thickness and the natural frequency of the laminate plate with clamped boundary conditions is greater than the simply supported plate as expected. The influence of circular cutout radius and temperature on the natural frequency is shown in \fref{fig:tempeffect}. It is seen that increasing the temperature, decreases the fundamental frequency while increasing the cutout radius, the fundamental frequency increases. The effect of the geometry of the cutout $d/e$ and the orientation of the cutout $\psi$ is shown in \fref{fig:elliporientation}. It can be seen that the orientation $\psi$ of the cutout and the size of the cutout $d/e$ has strong influence on the fundamental frequencies. The fundamental frequency decreases with increasing cutout size irrespective of the orientation of the cutout. With increasing orientation from 0$^\circ$ to 90$^\circ$, the frequency decreases and reaches minimum when the cutout is oriented at 60$^\circ$ and with further increase in the orientation, the frequency increases. Next, the effect of the cutout orientation and ply orientation on the fundamental frequency is studied. In this case, only one lamina is considered and the results are depicted in \fref{fig:plyorientation}. It is seen that for the ply orientation $\theta=$ 45$^\circ$, the frequency is symmetric with respect to cutout orientation $\psi=$ 45$^\circ$. 

\begin{figure}
\centering
\setlength\figureheight{10cm} 
\setlength\figurewidth{10cm}
\input{cvsfreq.tikz}
\caption{Normalized fundamental frequency as a function of moisture concentration, $C$ (\%) for a simply supported square laminated plate with (0$^\circ$/90$^\circ$/90$^\circ$/0$^\circ$) and a centrally located circular cutout, $r_o/a=$ 0.2 and for various plate thickness.}
\label{fig:moistureeffect}
\end{figure}
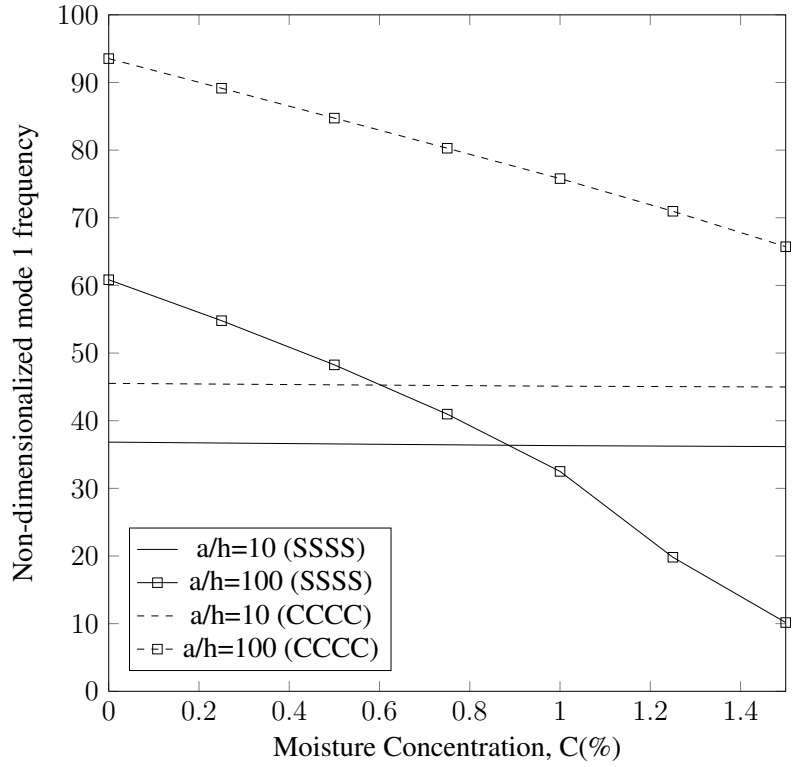

%% file: abVsfreq.tikz
%
%
%
%
\begin{tikzpicture}[scale=0.9]

\begin{axis}[%
width=\figurewidth,
height=\figureheight,
scale only axis,
xmin=0.5,
xmax=2.5,
xlabel={a/b},
ymin=10,
ymax=100,
ylabel={Non-dimensionalized mode 1 frequency},
legend pos = south east
]
\addplot [
color=black,
solid,
mark=o,
mark options={solid}
]
table[row sep=crcr]{
0.5 12.8372\\
0.6 14.5633\\
0.7 19.6506\\
0.8 24.1177\\
0.9 28.4312\\
1 36.3158\\
1.1 38.5304\\
1.2 43.583\\
1.3 46.4529\\
1.4 47.8899\\
1.5 49.6196\\
1.6 50.9856\\
1.7 51.794\\
1.8 52.8399\\
1.9 53.3035\\
2 55.3522\\
2.1 56.4842\\
2.2 56.8465\\
2.3 57.3098\\
2.4 58.2295\\
2.5 60.3249\\
};
\addlegendentry{a/h=10};

\addplot [
color=black,
dashed,
mark=square,
mark options={solid}
]
table[row sep=crcr]{
0.5 13.6108\\
0.6 15.9137\\
0.7 22.4877\\
0.8 29.1451\\
0.9 35.9858\\
1 48.5264\\
1.1 53.1757\\
1.2 63.4281\\
1.3 69.7361\\
1.4 72.1547\\
1.5 77.3701\\
1.6 80.4693\\
1.7 82.5274\\
1.8 84.5505\\
1.9 85.3157\\
2 88.6849\\
2.1 90.6357\\
2.2 90.8074\\
2.3 91.1821\\
2.4 92.4101\\
2.5 95.7119\\
};
\addlegendentry{a/h=20};

\end{axis}
\end{tikzpicture}%

%% file: cvsfreq.tikz
%
%
%
%
\begin{tikzpicture}[scale=0.9]

\begin{axis}[%
width=\figurewidth,
height=\figureheight,
scale only axis,
xmin=0,
xmax=1.5,
xlabel={Moisture Concentration, C(\%)},
ymin=0,
ymax=100,
ylabel={Non-dimensionalized mode 1 frequency},
legend pos= south west
]
\addplot [
color=black,
solid
]
table[row sep=crcr]{
0 36.8359\\
0.25 36.7011\\
0.5 36.5695\\
0.75 36.4411\\
1 36.3158\\
1.25 36.2443\\
1.5 36.1726\\
};
\addlegendentry{a/h=10 (SSSS)};

\addplot [
color=black,
solid,
mark=square,
mark options={solid}
]
table[row sep=crcr]{
0 60.8412\\
0.25 54.79\\
0.5 48.2572\\
0.75 40.9743\\
1 32.5061\\
1.25 19.8144\\
1.5 10.1659\\
};
\addlegendentry{a/h=100 (SSSS)};

\addplot [
color=black,
dashed
]
table[row sep=crcr]{
0 45.5191\\
0.25 45.4112\\
0.5 45.306\\
0.75 45.2033\\
1 45.1033\\
1.25 45.0449\\
1.5 44.9865\\
};
\addlegendentry{a/h=10 (CCCC)};

\addplot [
color=black,
dashed,
mark=square,
mark options={solid}
]
table[row sep=crcr]{
0 93.5338\\
0.25 89.1511\\
0.5 84.7343\\
0.75 80.2818\\
1 75.7913\\
1.25 70.9628\\
1.5 65.7252\\
};
\addlegendentry{a/h=100 (CCCC)};

\end{axis}
\end{tikzpicture}%

%% file: ravsfreq.tikz
%
%
%
%
\begin{tikzpicture}

\begin{axis}[%
width=\figurewidth,
height=\figureheight,
scale only axis,
xmin=0,
xmax=0.4,
xlabel={$\text{Cutout radius }{\text{r}_\text{o}}\text{/a}$},
ymin=10,
ymax=80,
ylabel={Non-dimensionalized mode 1 frequency},
legend pos = north west
]
\addplot [
color=black,
solid,
mark=square,
mark options={solid}
]
table[row sep=crcr]{
0 10.7978\\
0.1 23.29\\
0.2 36.8359\\
0.3 50.6925\\
0.4 76.6555\\
};
\addlegendentry{T=300K};

\addplot [
color=black,
solid,
mark=o,
mark options={solid}
]
table[row sep=crcr]{
0 10.1689\\
0.1 21.3265\\
0.2 33.6003\\
0.3 45.9787\\
0.4 69.1933\\
};
\addlegendentry{T=400K};

\end{axis}
\end{tikzpicture}%

%% file: devsfreq.tikz
%
%
%
%
\begin{tikzpicture}

\begin{axis}[%
width=\figurewidth,
height=\figureheight,
scale only axis,
xmin=0.4,
xmax=2,
xlabel={Size of the ellipse, d/e},
ymin=25,
ymax=55,
ylabel={Non-dimensionalized mode 1 frequency},
legend style={draw=black,fill=white,legend cell align=left}
]
\addplot [
color=black,
solid,
mark=square,
mark options={solid}
]
table[row sep=crcr]{
0.4 51.006\\
0.6 46.2135\\
0.8 40.8265\\
1 36.8359\\
1.2 32.7635\\
1.4 30.2332\\
1.6 29.7353\\
1.8 27.9236\\
2 27.9825\\
};
\addlegendentry{$\psi\text{=0}^\text{o}$};

\addplot [
color=black,
solid,
mark=o,
mark options={solid}
]
table[row sep=crcr]{
0.4 40.7115\\
0.6 40.7248\\
0.8 37.1946\\
1 36.8359\\
1.2 31.6415\\
1.4 31.4024\\
1.6 30.8635\\
1.8 30.3178\\
2 28.5979\\
};
\addlegendentry{$\psi\text{=45}^\text{o}$};

\end{axis}
\end{tikzpicture}%

%% file: deOrienvsfreq.tikz
%
%
%
%
\begin{tikzpicture}

\begin{axis}[%
width=\figurewidth,
height=\figureheight,
scale only axis,
xmin=-10,
xmax=90,
xlabel={$\text{Orientation of the ellipse, }\psi$},
ymin=38,
ymax=52,
ylabel={Non-dimensionalized mode 1 frequency},
legend style={draw=black,fill=white,legend cell align=left}
]
\addplot [
color=black,
solid,
mark=square,
mark options={solid}
]
table[row sep=crcr]{
0 51.006\\
15 44.7465\\
30 42.1063\\
45 40.7115\\
60 39.8132\\
75 40.4889\\
90 43.52\\
};
\addlegendentry{a/h=10, C=0\%, T=300K};

\end{axis}
\end{tikzpicture}%

%% file: singlevsfreq.tikz
%
%
%
%
\begin{tikzpicture}

\begin{axis}[%
width=\figurewidth,
height=\figureheight,
scale only axis,
xmin=0,
xmax=90,
xlabel={$\text{Orientation of the ellipse, }\psi$},
ymin=32,
ymax=46,
ylabel={Non-dimensionalized mode 1 frequency},
legend columns = -1,
legend to name = named,
]
\addplot [
color=black,
solid,
mark=square,
mark options={solid}
]
table[row sep=crcr]{
0 43.6729\\
15 38.4403\\
30 36.1639\\
45 34.8199\\
60 33.8174\\
75 34.1154\\
90 36.3363\\
};
\addlegendentry{$\theta\text{=0}^\text{o}$};

\addplot [
color=black,
solid,
mark=o,
mark options={solid}
]
table[row sep=crcr]{
0 45.9107\\
15 39.5857\\
30 36.9497\\
45 36.4158\\
60 36.9497\\
75 39.5857\\
90 45.9107\\
};
\addlegendentry{$\theta\text{=45}^\text{o}$};

\addplot [
color=black,
solid,
mark=diamond,
mark options={solid}
]
table[row sep=crcr]{
0 36.3363\\
15 34.1154\\
30 33.8174\\
45 34.8199\\
60 36.1639\\
75 38.4403\\
90 43.6729\\
};
\addlegendentry{$\theta\text{=90}^\text{o}$};

\end{axis}
\end{tikzpicture}%

\ref{named}

%% file: buck.tex
\paragraph{Buckling}

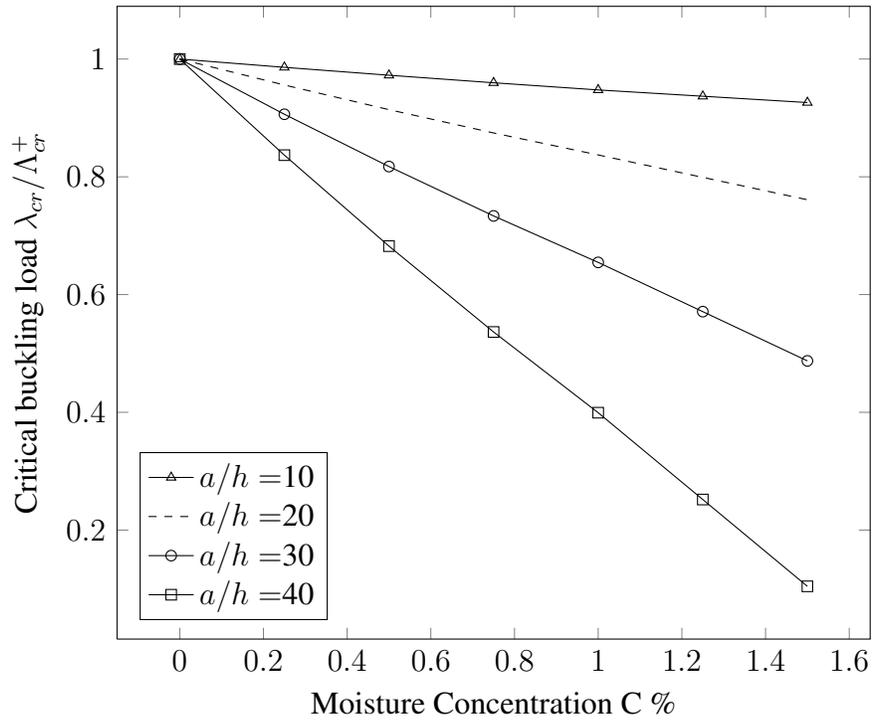
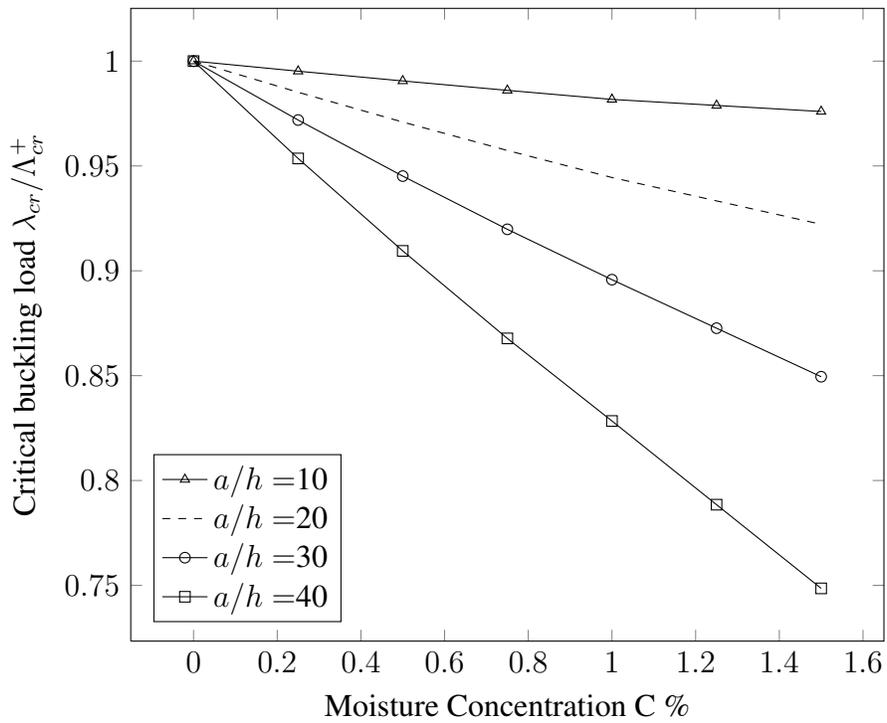
\begin{figure}
\centering
\subfigure[$a/b=$1]{
\pgfplotsset{compat=1.5}
\begin{tikzpicture}[scale=1.0]
\begin{axis}[
height=10cm,
xlabel={Moisture Concentration C \%},
ylabel={Critical buckling load $\lambda_{cr}/\Lambda_{cr}^+$},
legend pos= south west
]

\pgfplotstableread{abNoMoistNoCutout.dat}{\test}
\addplot [color=black,mark=triangle] table [x={C}, y={ah10}] {\test};
\addplot [dashed,black] table [x={C}, y={ah20}] {\test};
\addplot [color=black,mark=o] table [x={C}, y={ah30}] {\test};
\addplot [color=black,mark=square] table [x={C}, y={ah40}] {\test};
\legend{$a/h=$10, $a/h=$20,$a/h=$30,$a/h=$40} 
\end{axis}
\end{tikzpicture}}

\subfigure[$a/b=$2]{
\pgfplotsset{compat=1.5}
\begin{tikzpicture}[scale=1.0]
\begin{axis}[
height=10cm,
xlabel={Moisture Concentration C \%},
ylabel={Critical buckling load $\lambda_{cr}/\Lambda_{cr}^+$},
legend pos= south west
]

\pgfplotstableread{ab2NoMoistNoCutout.dat}{\test}
\addplot [color=black,mark=triangle] table [x={C}, y={ah10}] {\test};
\addplot [dashed,black] table [x={C}, y={ah20}] {\test};
\addplot [color=black,mark=o] table [x={C}, y={ah30}] {\test};
\addplot [color=black,mark=square] table [x={C}, y={ah40}] {\test};
\legend{$a/h=$10, $a/h=$20,$a/h=$30,$a/h=$40} 
\end{axis}
\end{tikzpicture}}
\caption{Influence of moisture concentration (C \%) on the normalized critical buckling load $\lambda_{cr}/\Lambda_{cr}^+$ for a square and a rectangular laminated plate (0$^\circ$/90$^\circ$/90$^\circ$/0$^\circ$). $\Lambda_{cr}^+$ is the critical buckling load of a composite laminate without a cutout.}
\label{fig:mositeffectbuck}
\end{figure}

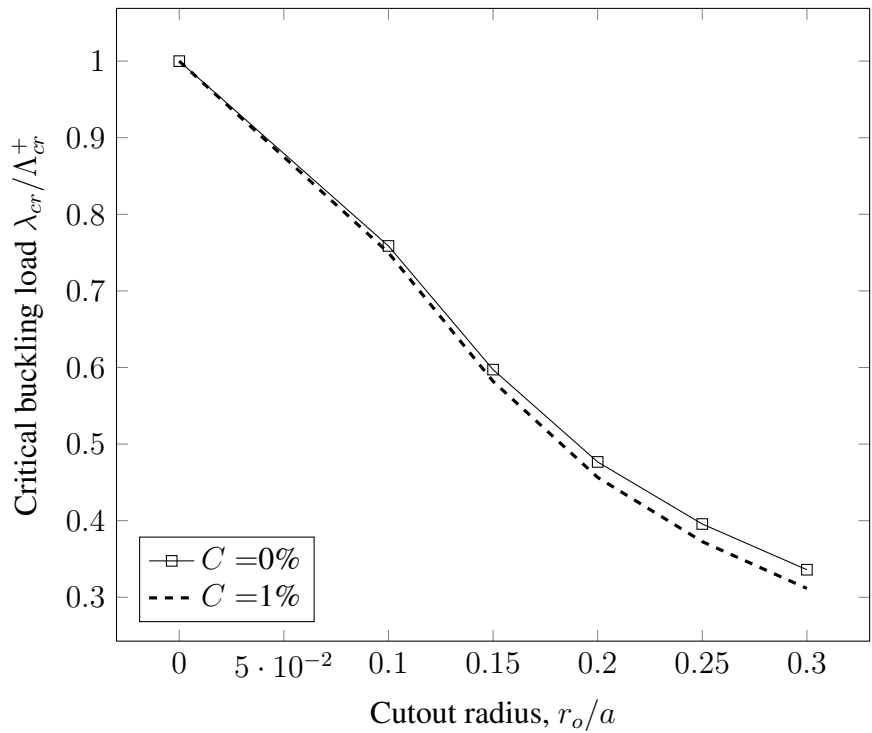
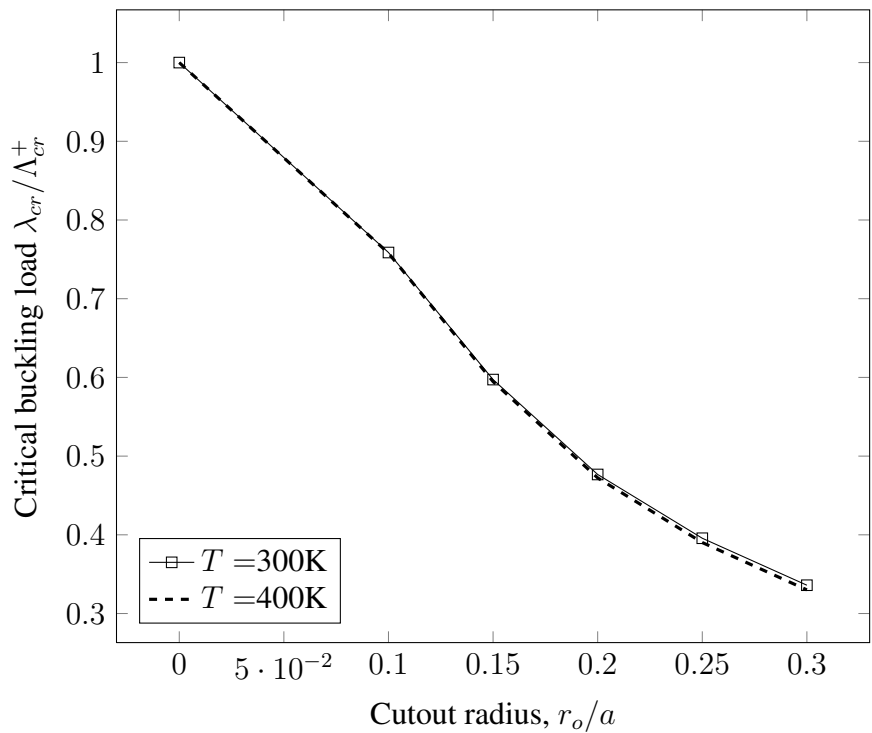
\begin{figure}
\centering
\subfigure[Influence of Moisture Concentration, $C$]{
\pgfplotsset{compat=1.5}
\begin{tikzpicture}[scale=1.0]
\begin{axis}[
height=10cm,
xlabel={Cutout radius, $r_o/a$},
ylabel={Critical buckling load $\lambda_{cr}/\Lambda_{cr}^+$},
legend pos= south west
]

\pgfplotstableread{raMoistCutout.dat}{\test}
\addplot [color=black,mark=square] table [x={ra}, y={c0}] {\test};
\addplot [dashed,black,very thick] table [x={ra}, y={c1}] {\test};
\legend{$C=$0\%, $C=$1\%} 
\end{axis}
\end{tikzpicture}}

\subfigure[Influence of Temperature, $T$]{
\pgfplotsset{compat=1.5}
\begin{tikzpicture}[scale=1.0]
\begin{axis}[
height=10cm,
xlabel={Cutout radius, $r_o/a$},
ylabel={Critical buckling load $\lambda_{cr}/\Lambda_{cr}^+$},
legend pos= south west
]

\pgfplotstableread{raTempCutout.dat}{\test}
\addplot [color=black,mark=square] table [x={ra}, y={T300}] {\test};
\addplot [dashed,black,very thick] table [x={ra}, y={T400}] {\test};
\legend{$T=$300K, $T=$400K} 
\end{axis}
\end{tikzpicture}}
\caption{Influence of moisture concentration (C \%), Temperature and a centrally located circular cutout on the normalized critical buckling load $\lambda_{cr}/\Lambda_{cr}^+$ for a square laminated plate with (0$^\circ$/90$^\circ$/90$^\circ$/0$^\circ$). $\Lambda_{cr}^+$ is the critical buckling load of a composite laminate without a cutout.}
\label{fig:cutoutradimoistempbuck}
\end{figure}

\begin{figure}
\centering
\pgfplotsset{compat=1.5}
\begin{tikzpicture}[scale=1.0]
\begin{axis}[
height=10cm,
xlabel={Cutout radius, $r_o/a$},
ylabel={Critical buckling load $\lambda_{cr}/\Lambda_{cr}^+$},
legend pos= south west
]
\pgfplotstableread{BcEffectT325.dat}{\test}
\addplot [color=black,mark=square] table [x={ra}, y={SSSS}] {\test};
\addplot [dashed,black] table [x={ra}, y={CCCC}] {\test};
\legend{SSSS, CCCC} 
\end{axis}
\end{tikzpicture}
\caption{Effect of various boundary conditions and the size of a centrally located circular cutout on the normalized critical buckling load $\lambda_{cr}/\Lambda_{cr}^+$ for a square composite laminate with (0$^\circ$/90$^\circ$/90$^\circ$/0$^\circ$). $\Lambda_{cr}^+$ is the critical buckling load of a composite laminate without a cutout.}
\label{fig:effectbcs}
\end{figure}
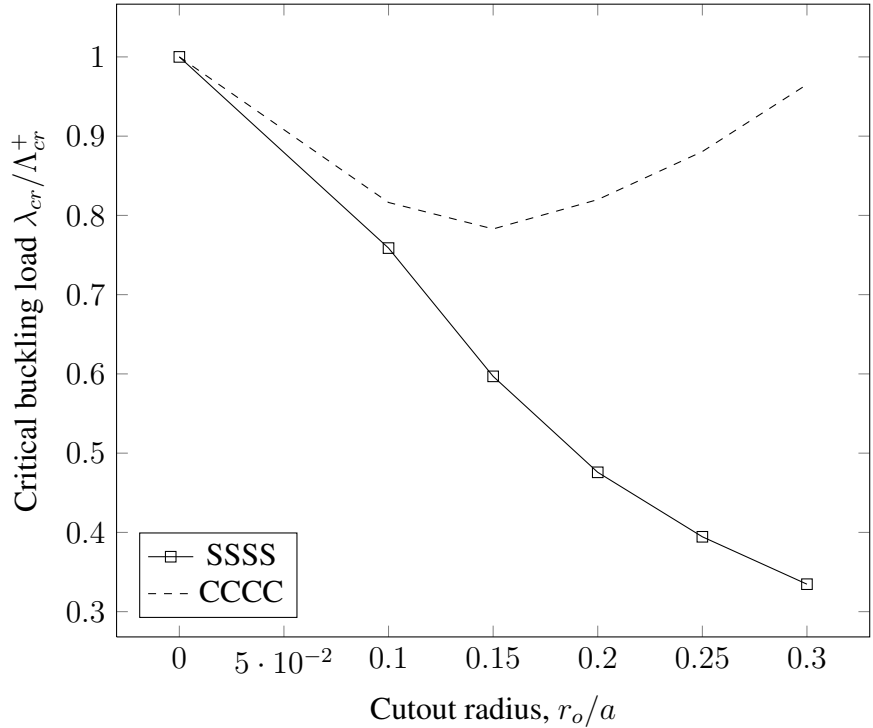

Next, the effect of moisture concentration, the cutout and temperature on the critical load is numerically investigated. A laminated plate with ply sequence (0$^\circ$/90$^\circ$/90$^\circ$/0$^\circ$) is considered for the analysis. In all cases, we present the non dimensionalized critical load as, unless specified otherwise:
\begin{equation}
\overline{N}_{xx}^{\rm cr} = \lambda_{cr}/\Lambda_{cr}^+
\label{eqn:nondimbuck}
\end{equation}
where $N_{xx}^{\ast {\rm cr}}$ is the critical load of the laminate plate without a cutout and with moisture concentration $C=$ 0\% and temperature $T=$ 300K. \fref{fig:mositeffectbuck} shows the effect of moisture concentration, the plate thickness and the aspect ratio on the critical buckling load for a simply supported laminated plate. It is seen that with increasing moisture concentration and $a/h$, the critical buckling load decreases, whilst with increasing plate aspect ratio, the critical load increases. The effect of cutout radius $r_o/a$, the moisture concentration $C$ and the thermal gradient on the critical buckling load is shown in \fref{fig:cutoutradimoistempbuck}. It is seen that with increasing the cutout radius, the moisture concentration and thermal gradient, the critical buckling load decreases. This can be attributed to the stiffness degradation. \fref{fig:effectbcs} shows the influence of various boundary conditions and the size of a centrally located circular cutout on the normalized critical buckling load $\lambda_{cr}/\Lambda_{cr}^+$ for a square composite laminate with (0$^\circ$/90$^\circ$/90$^\circ$/0$^\circ$). It is seen that with increasing cutout radius, the critical buckling load decreases for a simply supported plate, while for a clamped plate, the critical buckling load first decreases and then increases.

